\documentclass [12pt]{article}
\usepackage{amssymb}
\def\lanbox{\hbox{$\, \vrule height 0.25cm width 0.25cm depth 0.01cm \,$}}

\usepackage{color}

\begin{document}

\centerline{\Large\bf Existence of solutions for some systems of superdiffusive}
\centerline{\Large\bf integro-differential equations in
population dynamics}
\centerline{\Large\bf depending on the natality and mortality rates}

\bigskip

\centerline{Vitali Vougalter}

\bigskip

\centerline{Department of Mathematics, University
of Toronto}

\centerline{Toronto, Ontario, M5S 2E4, Canada}

\centerline{ e-mail: vitali@math.toronto.edu}

\bigskip
\bigskip
\bigskip

\noindent {\bf Abstract.}
We prove the existence of stationary solutions for some systems of
reaction-diffusion type equations with superdiffusion in the corresponding
$H^{2}$ spaces.
Our method is based on the fixed point theorem when the elliptic problems
contain first order differential operators with and without the Fredholm
property, which may depend on the outcome of the competition between the
natality and the mortality rates contained in the equations of our systems.

\bigskip
\bigskip



\bigskip
\bigskip
\bigskip
\bigskip

\section{Introduction}

\bigskip

\noindent
We recall that a linear operator $L$, which acts from a Banach space $E$
into another Banach space $F$ has the Fredholm property if its
image is closed, the dimension of its kernel and the codimension of its
image are finite. Consequently, the equation $Lu=f$ is solvable if and 
only if $\phi_{k}(f)=0$ for a finite number of functionals $\phi_{k}$ from
the dual space $F^{*}$. Such properties of the Fredholm operators are widely
used in different methods of the linear and nonlinear analysis.

\noindent
Elliptic problems studied in bounded domains with a sufficiently smooth
boundary satisfy the Fredholm property when the ellipticity condition,
proper ellipticity and Lopatinskii conditions are fulfilled (see e.g. 
~\cite{A97}, ~\cite{LM68}, ~\cite{V68}). This is the main result of the
theory of linear elliptic equations. When dealing with unbounded domains,
such conditions may not be sufficient and the Fredholm property may not
be satisfied. For instance, for the Laplace operator, $Lu=\Delta u$ considered
in ${\mathbb R}^{d}$ the Fredholm property fails to hold when the problem is
studied either in H\"older spaces, so that $L: {\mathbb C}^{2+\alpha}
({\mathbb R}^{d})\to {\mathbb C}^{\alpha}({\mathbb R}^{d})$ or in Sobolev
spaces, $L: H^{2}({\mathbb R}^{d})\to L^{2}({\mathbb R}^{d})$.

\noindent
For linear elliptic equations discussed in unbounded domains the Fredholm 
property is satisfied if and only if, additionally to the conditions given
above, the limiting operators are invertible (see ~\cite{V11}). In certain
simple cases, the limiting operators can be constructed explicitly. For
instance, when
$$
Lu=a(x)u''+b(x)u'+c(x)u, \quad x\in {\mathbb R},
$$
with the coefficients of the operator having limits at infinity,
$$
a_{\pm}=\hbox{lim}_{x\to \pm \infty}a(x), \quad 
b_{\pm}=\hbox{lim}_{x\to \pm \infty}b(x), \quad
c_{\pm}=\hbox{lim}_{x\to \pm \infty}c(x),
$$
the limiting operators are equal to
$$
L_{\pm}u=a_{\pm}u''+b_{\pm}u'+c_{\pm}u.
$$
Because the coefficients here are constants, the essential
spectrum of the operator, which is the set of complex numbers $\lambda$
for which the operator $L-\lambda$ does not possess the Fredholm property,
can be found explicitly using the standard Fourier transform, so that
$$
\lambda_{\pm}(\xi)=-a_{\pm}\xi^{2}+b_{\pm}i\xi+c_{\pm}, \quad \xi\in {\mathbb R}.
$$
The limiting operators are invertible if and only if the origin is not
contained in the essential spectrum. 

\noindent
For the general elliptic problems the analogous assertions hold true. The
Fredholm property is satisfied when the origin does not belong to the essential
spectrum or if the limiting operators are invertible. These conditions may not
be written explicitly.

\noindent
For the non-Fredholm operators we may not apply the standard solvability 
relations and in a general case the solvability conditions are unknown. 
However, the solvability relations were derived recently for certain classes
of operators. For instance, consider the following equation
\begin{equation}
\label{eq1}
Lu\equiv \Delta u+au=f
\end{equation}
in ${\mathbb R}^{d}, \ d\in {\mathbb N}$ with a positive constant $a$. 
Here the operator $L$ and its limiting operators coincide. The corresponding 
homogeneous problem has a nontrivial bounded solution, so that the 
Fredholm property is not satisfied.
Since the differential operator contained in (\ref{eq1}) has 
constant coefficients, we are able to find the solution explicitly by means of
the standard Fourier transform. In Lemmas 5 and 6 of ~\cite{VV113} we obtained
the following solvability conditions. Let $f(x)\in L^{2}({\mathbb R}^{d})$ and 
$xf(x)\in L^{1}({\mathbb R}^{d})$. Then equation (\ref{eq1}) admits a unique 
solution in $H^{2}({\mathbb R}^{d})$ if and only if
$$
\Bigg(f(x), \frac{e^{ipx}}{(2\pi)^{\frac{d}{2}}}\Bigg)_{L^{2}({\mathbb R}^{d})}=0, \quad
p\in S_{\sqrt{a}}^{d} \quad a.e.
$$
Here and below $S_{r}^{d}$ stands for the sphere in ${\mathbb R}^{d}$
of radius $r$ centered at the origin. Hence, although the Fredholm property
is not satisfied for our operator, we are able to formulate the solvability
relations similarly. Note that such similarity is only formal
since the range of the operator is not closed. In the situation
when the operator contains a scalar potential, so that
$$
Hu \equiv \Delta u+b(x)u=f,
$$
we cannot use the standard Fourier transform directly. However, the
solvability conditions in three dimensions can be derived by means of the
spectral and the scattering theory of Schr\"odinger type operators
(see ~\cite{VV08}).
Similarly to the constant coefficient case, the solvability relations are
expressed
in terms of orthogonality to the solutions of the adjoint homogeneous problem.
The solvability conditions for several other examples of non-Fredholm
linear elliptic operators were obtained (see ~\cite{EV21}, ~\cite{EV22},
~\cite{V11}, ~\cite{VKMP02}, ~\cite{VV101}, ~\cite{VV08}, ~\cite{VV113}).

\noindent
Solvability relations are important in the analysis of the nonlinear elliptic
equations. When the non-Fredholm operators are involved, in spite of a certain
progress in the studies of the linear equations, the nonlinear non-Fredholm
operators were analyzed only in few examples (see ~\cite{DMV05}, ~\cite{DMV08},
~\cite{DMV08b}, ~\cite{EV211}, ~\cite{EV22}, ~\cite{EV221},
~\cite{VV113}).
Clearly, such situation can be explained by the fact that the majority of 
the methods of linear and nonlinear analysis are based on the Fredholm property.
Fredholm structures, topological invariants and applications were discussed
in ~\cite{E09}. Symmetrization and stabilization of solutions of nonlinear
elliptic equations were covered in ~\cite{E18}.

\noindent
The present article deals with the  certain systems of nonlinear integro-
differential reaction-diffusion type equations, for 
which the Fredholm property may not be satisfied:
\begin{equation}
\label{id1}
\frac{\partial u_{k}}{\partial t}=-\sqrt{-\Delta}u_{k}+\int_{\Omega}G_{k}(x-y)
F_{k}(u_{1}(y,t), u_{2}(y,t),..., u_{N_{2}}(y,t), y)dy
+a_{k}u_{k}
\end{equation}
with $1\leq k\leq N_{1}$ and
\begin{equation}
\label{id11}
\frac{\partial u_{k}}{\partial t}=-\sqrt{-\Delta}u_{k}+\int_{\Omega}G_{k}(x-y)
F_{k}(u_{1}(y,t), u_{2}(y,t),..., u_{N_{2}}(y,t), y)dy
-a_{k}u_{k}
\end{equation}
with $N_{1}+1\leq k\leq N_{2}$. In system (\ref{id1}), (\ref{id11}) the constants
$\{a_{k}\}_{k=1}^{N_{2}}$ are nonnegative when $1\leq k\leq N_{1}$ and they are
strictly positive for $N_{1}+1\leq k\leq N_{2}$. The domain here
$\Omega \subseteq {\mathbb R}^{d}, \quad 1\leq d\leq 3$. These are the more
physically relevant dimensions. Let us note that equations (\ref{id1})
describe the cases in the Population Dynamics in the Mathematical Biology when
the natality rates are higher than the mortality ones for
$a_{k}>0$ and the situations when the mortality and natality rates balance each
other if $a_{k}=0$. On the other hand, equations (\ref{id11}) are crucial
for understanding the cases when the mortality rates are higher than the
natality ones.

\noindent
The operator $\sqrt{-\Delta}$ in the equations of the system above  is defined
by means of the spectral calculus and is actively used, for example in the
studies of the superdiffusion problems (see e.g. ~\cite{MK00}, ~\cite{VV18}
and the references therein). Superdiffusion can be described as a random
process of the particle motion characterized by the probability density
distribution of the jump length. The moments of this density distribution
are finite in the case of the normal diffusion, but this is not the case for
superdiffusion. The asymptotic behavior at the infinity of the probability
density function determines the value of the power of the negative Laplacian
(see ~\cite{MK00}).

\noindent
In the Population Dynamics the
integro-differential problems are used to describe the biological 
systems with the intra-specific competition and the nonlocal consumption 
of resources (see e.g. \cite{ABVV10}, \cite{BNPR09}, \cite{GVA06}, 
\cite{VV112}). The stability issues for the travelling fronts of reaction-
diffusion type equations with the essential spectrum of the linearized operator
crossing the imaginary axis were also discussed in ~\cite{BGS09} and 
~\cite{GS07}. The single equation analogous to (\ref{id1}) was considered
in \cite{VV18}. The work ~\cite{SZ10} is devoted to the spreading speeds
for monostable equations with nonlocal dispersal in space periodic habitats.
A discretization method for nonlocal diffusion type equations was developed
in ~\cite{MORV22}. The article ~\cite{KKM16} deals with the
weak error for continuous time Markov chains related to fractional
in time P(I)DEs.

\noindent
Let us introduce
$$
{\cal F}(u,x):=(F_{1}(u, x), F_{2}(u, x), ..., F_{N_{2}}(u, x))^{T}.
$$
The nonlinear terms of  system (\ref{id1}), (\ref{id11}) will satisfy the
following regularity requirements.

\bigskip

\noindent
{\bf Assumption 1.1.} {\it Let $1\leq k\leq N_{2}$. Functions
$F_{k}(u,x): {\mathbb R}^{N_{2}}\times\Omega \to {\mathbb R}$ are satisfying the
Caratheodory condition (see ~\cite{K64}), such that
\begin{equation}
\label{ub1}
|{\cal F}(u, x)|_{{\mathbb R}^{N_{2}}}\leq K|u|_{{\mathbb R}^{N_{2}}}+h(x) \quad for 
\quad u\in {\mathbb R}^{N_{2}}, \quad x\in \Omega,
\end{equation}
with a constant $K>0$ and $h(x):\Omega\to {\mathbb R}^{+}, \quad
h(x)\in L^{2}(\Omega)$. Moreover, they are Lipschitz continuous functions, 
such that for any  $u^{(1),(2)}\in{\mathbb R}^{N_{2}}, \quad x\in \Omega:$
\begin{equation}
\label{lk1}
|{\cal F}(u^{(1)}, x)-{\cal F}(u^{(2)}, x)|_{{\mathbb R}^{N_{2}}}\leq L |u^{(1)}-u^{(2)}|_
{{\mathbb R}^{N_{2}}},
\end{equation}
with a constant $L>0$.}

\bigskip

\noindent
Here and further down we use the notations for a vector
$u:=(u_{1},u_{2},...,u_{N_{2}})^{T}\in {\mathbb R}^{N_{2}}$ and its norm
$|u|_{{\mathbb R}^{N_{2}}}:=\sqrt{\sum_{k=1}^{N_{2}}u_{k}^{2}}$.

\noindent
The solvability of a
local elliptic problem in a bounded domain in ${\mathbb R}^{N}$ was considered
in ~\cite{BO86}. The nonlinear function involved there was allowed to have a
sublinear growth.

\noindent
Obviously, the stationary solutions of system (\ref{id1}), (\ref{id11}), if any
exist, will satisfy the following system of nonlocal elliptic equations
$$
-\sqrt{-\Delta}u_{k}+\int_{\Omega}G_{k}(x-y)F_{k}(u_{1}(y),u_{2}(y),...,u_{N_{2}}(y),
y)dy+a_{k}u_{k}=0, \quad a_{k}\geq 0,
$$
for $1\leq k\leq N_{1}$, 
$$
-\sqrt{-\Delta}u_{k}+\int_{\Omega}G_{k}(x-y)F_{k}(u_{1}(y),u_{2}(y),...,u_{N_{2}}(y),
y)dy-a_{k}u_{k}=0, \quad a_{k}>0
$$
for $N_{1}+1\leq k\leq N_{2}$.

\noindent
Let us introduce the auxiliary semi-linear problem for the technical purposes,
namely
\begin{equation}
\label{ae1}
\sqrt{-\Delta}u_{k}-a_{k}u_{k}=\int_{\Omega}G_{k}(x-y)F_{k}(v_{1}(y),v_{2}(y),...,
v_{N_{2}}(y), y)dy, \quad a_{k}\geq 0
\end{equation}
if $1\leq k\leq N_{1}$,
\begin{equation}
\label{ae2}
\sqrt{-\Delta}u_{k}+a_{k}u_{k}=\int_{\Omega}G_{k}(x-y)F_{k}(v_{1}(y),v_{2}(y),...,
v_{N_{2}}(y), y)dy, \quad a_{k}>0
\end{equation}
if $N_{1}+1\leq k\leq N_{2}$.

\noindent
We denote
\begin{equation}
\label{ip}  
(f_{1}(x),f_{2}(x))_{L^{2}(\Omega)}:=\int_{\Omega}f_{1}(x)\bar{f_{2}}(x)dx,
\end{equation}
with a slight abuse of notations in the situation when such functions do not
belong to $L^{2}(\Omega)$, like for instance those used in the 
orthogonality conditions of the assumption below. Indeed, if 
$f_{1}(x)\in L^{1}(\Omega)$ and $f_{2}(x)$ is bounded in $\Omega$, then the
integral in the right side of (\ref{ip}) makes sense.

\noindent
Let us begin the article by considering the whole space case, so that
$\Omega={\mathbb R}^{d}$ and the corresponding Sobolev space is
equipped with the norm
$$
\|u\|_{H^{2}({\mathbb R}^{d}, \ {\mathbb R}^{N_{2}})}^{2}:=\sum_{k=1}^{N_{2}}\|u_{k}\|_{H^{2}
({\mathbb R}^{d})}^{2}=\sum_{k=1}^{N_{2}}\{\|u_{k}\|_{L^{2}({\mathbb R}^{d})}^{2}+
\|\Delta u_{k}\|_{L^{2}({\mathbb R}^{d})}^{2}\},
$$
where $u(x): {\mathbb R}^{d}\to {\mathbb R}^{N_{2}}$.

\noindent
The primary obstacle in solving system (\ref{ae1}), (\ref{ae2}) is that
operators
$\sqrt{-\Delta}-a_{k}: H^{1}({\mathbb R}^{d})\to L^{2}({\mathbb R}^{d}), \
a_{k}\geq 0$ contained in its first part do not satisfy the Fredholm property.

\noindent
The analogous situations in linear problems, which can be self-adjoint or 
non self-adjoint containing the non-Fredholm second and higher order 
differential operators or even systems of equations involving the non-Fredholm 
operators have been studied actively in recent years (see ~\cite{EV21},
~\cite{EV22}, ~\cite{V11}, ~\cite{VKMP02}, ~\cite{VV101}, ~\cite{VV08},
~\cite{VV113}).

\noindent
We manage to demonstrate that system (\ref{ae1}), (\ref{ae2}) defines
a map $T_{a}: H^{2}({\mathbb R}^{d}, \ {\mathbb R}^{N_{2}})\to
H^{2}({\mathbb R}^{d}, \ {\mathbb R}^{N_{2}})$,
which is a strict contraction under the given technical assumptions. It
gives a solution of the considered problem. The fact that this map is well
defined is established in the proof of Theorem 1.3 further down.

\noindent
Let us formulate the following conditions on the integral kernels contained in
the nonlocal parts of the system of equations (\ref{ae1}), (\ref{ae2}).

\bigskip

\noindent
{\bf Assumption 1.2.} {\it  Let $G_{k}(x): {\mathbb R}^{d}\to {\mathbb R}, \ 
G_{k}(x)\in W^{1, 1}({\mathbb R}^{d}), \ 1\leq k\leq N_{2}, \ 1\leq d\leq 3$ and
$1\leq m\leq N_{1}-1, \ m\in {\mathbb N}$ with $N_{1}\geq 2, \ N_{2}>N_{1}$.

\medskip

\noindent
I) Let $a_{k}>0, \ 1\leq k\leq m$, assume that 
$xG_{k}(x)\in L^{1}({\mathbb R}^{d})$ and 
\begin{equation}
\label{or1}
\Bigg(G_{k}(x), \frac{e^{\pm ia_{k}x}}{\sqrt{2\pi}}\Bigg)_{L^{2}({\mathbb R})}=0
\quad when \quad d=1.
\end{equation}
\begin{equation}
\label{or2}
\Bigg(G_{k}(x),\frac{e^{ipx}}{(2\pi)^{\frac{d}{2}}}\Bigg)_{L^{2}({\mathbb R}^{d})}=0 \quad
for \quad p\in S_{a_{k}}^{d} \quad for \quad d=2,3.
\end{equation}

\medskip

\noindent
II) Let $a_{k}=0, \ m+1\leq k\leq N_{1}$, assume that 
$xG_{k}(x)\in L^{1}({\mathbb R}^{d})$ and 
\begin{equation}
\label{or3}
(G_{k}(x), 1)_{L^{2}({\mathbb R}^{d})}=0. 
\end{equation}

\medskip

\noindent
III) Let $a_{k}>0, \ N_{1}+1\leq k\leq N_{2}$.
}

\bigskip

\noindent
We use the hat symbol here and further down to designate the standard 
Fourier transform, namely
\begin{equation}
\label{ft}  
\widehat{G_{k}}(p):=\frac{1}{(2\pi)^{\frac{d}{2}}}\int_{{\mathbb R}^{d}}G_{k}(x)
e^{-ipx}dx, \quad p\in {\mathbb R}^{d}.
\end{equation}
Hence,
\begin{equation}
\label{fub}  
\|\widehat{G_{k}}(p)\|_{L^{\infty}({\mathbb R}^{d})}\leq \frac{1}{(2\pi)^{\frac{d}{2}}}
\|G_{k}\|_{L^{1}({\mathbb R}^{d})}
\end{equation}
and
\begin{equation}
\label{fubg}  
\|p\widehat{G_{k}}(p)\|_{L^{\infty}({\mathbb R}^{d})}\leq \frac{1}{(2\pi)^{\frac{d}{2}}}
\|\nabla G_{k}\|_{L^{1}({\mathbb R}^{d})}.
\end{equation}
Let us define the following technical expressions
\begin{equation}
\label{Mka}
M_{k}:=\hbox{max}\Bigg\{\Bigg\| \frac{\widehat{G_{k}}(p)}{|p|-a_{k}}\Bigg\|_
{L^{\infty}({\mathbb R}^{d})}, \ \Bigg\| \frac{p^{2}\widehat{G_{k}}(p)}{|p|-a_{k}}
\Bigg\|_{L^{\infty}({\mathbb R}^{d})}\Bigg\}, \quad 1\leq k\leq m. 
\end{equation}
\begin{equation}
\label{Mk0}
M_{k}:=\hbox{max}\Bigg\{\Bigg\| \frac{\widehat{G_{k}}(p)}{|p|}\Bigg\|_
{L^{\infty}({\mathbb R}^{d})}, \ \Bigg\|p\widehat{G_{k}}(p)
\Bigg\|_{L^{\infty}({\mathbb R}^{d})}\Bigg\}, \quad m+1\leq k\leq N_{1}. 
\end{equation}
\begin{equation}
\label{Mka+}
M_{k}:=\hbox{max}\Bigg\{\Bigg\| \frac{\widehat{G_{k}}(p)}{|p|+a_{k}}\Bigg\|_
{L^{\infty}({\mathbb R}^{d})}, \ \Bigg\| \frac{p^{2}\widehat{G_{k}}(p)}{|p|+a_{k}}
\Bigg\|_{L^{\infty}({\mathbb R}^{d})}\Bigg\}, \quad N_{1}+1\leq k\leq N_{2}. 
\end{equation}
Clearly, quantities (\ref{Mka}) and (\ref{Mk0}) are finite by means of
Lemma A1 in one dimension and Lemma A2 when $d=2,3$ of \cite{VV18} under
Assumption 1.2 above.

\noindent
It can be trivially checked that (\ref{Mka+}) are finite
as well. Indeed, for $N_{1}+1\leq k\leq N_{2}$ by virtue of (\ref{fub}), we
obtain
$$
\Bigg|\frac{\widehat{G_{k}}(p)}{|p|+a_{k}}\Bigg|\leq
\frac{|\widehat{G_{k}}(p)|}{a_{k}}\leq
\frac{1}{(2\pi)^{\frac{d}{2}}a_{k}}\|G_{k}\|_{L^{1}({\mathbb R}^{d})}<\infty
$$
as we assume. Similarly, by means of (\ref{fubg})
$$
\Bigg|\frac{p^{2}\widehat{G_{k}}(p)}{|p|+a_{k}}\Bigg|\leq
|p\widehat{G_{k}}(p)|\leq
\frac{1}{(2\pi)^{\frac{d}{2}}}\|\nabla G_{k}\|_{L^{1}({\mathbb R}^{d})}<\infty
$$
due to the given condition. Therefore,
$M_{k}<\infty$ for $N_{1}+1\leq k\leq N_{2}$ as well.

\noindent
This allows us to introduce
\begin{equation}
\label{M}
M:=\hbox{max}M_{k}, \quad 1\leq k\leq N_{2},
\end{equation}
where $M_{k}$ are given by (\ref{Mka}), (\ref{Mk0}) and (\ref{Mka+}).

\noindent
Our first main result is as follows.

\bigskip

\noindent
{\bf Theorem 1.3.}  {\it Let $\Omega={\mathbb R}^{d}, \ d=1,2,3$,
Assumptions 1.1 and 1.2 hold and $\sqrt{2}(2\pi)^{\frac{d}{2}}ML<1$.
 
\medskip

\noindent
Then the map $T_{a}v=u$ on $H^{2}({\mathbb R}^{d}, \ {\mathbb R}^{N_{2}})$ 
defined by system (\ref{ae1}), (\ref{ae2}) has a unique
fixed point $v_{a}(x):{\mathbb R}^{d}\to {\mathbb R}^{N_{2}}$, which is the only
stationary solution of problem (\ref{id1}), (\ref{id11}) in 
$H^{2}({\mathbb R}^{d}, \ {\mathbb R}^{N_{2}})$.

\medskip

\noindent
This fixed point $v_{a}(x)$ is nontrivial provided the intersection of 
supports of the Fourier images of functions 
$supp\widehat{F_{k}(0,x)}(p)\cap supp \widehat{G_{k}}(p)$ is a set 
of nonzero Lebesgue measure in ${\mathbb R}^{d}$ for some
$1\leq k\leq N_{2}$.}

\bigskip

\noindent
We turn our attention to the studies of the analogous system of equations
on the interval $\Omega=I:=[0, \ 2\pi]$ with periodic boundary conditions
for the solution vector function and its first derivative. Let us assume the 
following about the integral kernels contained in the nonlocal parts of system 
(\ref{ae1}), (\ref{ae2}) in such case.

\bigskip

\noindent
{\bf Assumption 1.4.} {\it  Let $\displaystyle{G_{k}(x): I\to {\mathbb R}, \
G_{k}(x)\in C(I), \ \frac{dG_{k}(x)}{dx}\in L^{1}(I)}$ with

\noindent
$G_{k}(0)=G_{k}(2\pi), \ 1\leq k\leq N_{2}$, where
$N_{1}\geq 3, \ 1\leq m<q\leq N_{1}-1, \ m,q\in {\mathbb N}$ and
$N_{2}>N_{1}$.

\medskip

\noindent
I) Let $a_{k}>0$ and $a_{k}\neq n, \quad n\in {\mathbb N}$ \quad if \quad
$1\leq k\leq m$.

\medskip

\noindent
II) Let $a_{k}=n_{k}, \ n_{k}\in {\mathbb N}$ and
\begin{equation}
\label{or4}
\Bigg(G_{k}(x), \frac{e^{\pm in_{k}x}}{\sqrt{2\pi}}\Bigg)_{L^{2}(I)}=0 \quad if \quad
m+1\leq k \leq q.
\end{equation}

\medskip

\noindent
III) Let $a_{k}=0$ and
\begin{equation}
\label{or5}
(G_{k}(x), 1)_{L^{2}(I)}=0 \quad for \quad q+1\leq k\leq N_{1}.
\end{equation}

\medskip

\noindent
IV) Let $a_{k}>0, \ N_{1}+1\leq k\leq N_{2}$.

\medskip

\noindent
Let $F_{k}(u,0)=F_{k}(u,2\pi)$ for $u\in {\mathbb R}^{N_{2}}$ and
$1\leq k\leq N_{2}$.}

\bigskip

\noindent
For the function on the $[0, 2\pi]$ interval,
$G_{k}(x): I\to {\mathbb R}, \ G_{k}(0)=G_{k}(2\pi)$, we
use the Fourier transform given by
\begin{equation}
\label{fti}  
G_{k, \ n}:=\int_{0}^{2\pi}G_{k}(x)\frac{e^{-inx}}{\sqrt{2\pi}}dx, \quad
n\in {\mathbb Z},
\end{equation}
so that
$$
G_{k}(x)=\sum_{n=-\infty}^{\infty}G_{k, \ n}\frac{e^{inx}}{\sqrt{2\pi}}.
$$
Evidently, the upper bounds
\begin{equation}
\label{fiub}
\|G_{k, \ n}\|_{l^{\infty}}\leq \frac{1}{\sqrt{2\pi}}\|G_{k}\|_{L^{1}(I)}, \quad
\|nG_{k, \ n}\|_{l^{\infty}}\leq \frac{1}{\sqrt{2\pi}}\|G_{k}'(x)\|_{L^{1}(I)}
\end{equation}
hold.

\noindent
Let us introduce the following auxiliary quantities
\begin{equation}
\label{Pk1}
P_{k}:=\hbox{max}\Bigg\{ \Bigg\|\frac{G_{k, \ n}}{|n|-a_{k}}\Bigg\|_{l^{\infty}},
\Bigg\|\frac{n^{2}G_{k, \ n}}{|n|-a_{k}}\Bigg\|_{l^{\infty}} \Bigg \}, \quad 1\leq k
\leq m.
\end{equation}
\begin{equation}
\label{Pk2}
P_{k}:=\hbox{max}\Bigg\{ \Bigg\|\frac{G_{k, \ n}}{|n|-n_{k}}\Bigg\|_{l^{\infty}},
\Bigg\|\frac{n^{2}G_{k, \ n}}{|n|-n_{k}}\Bigg\|_{l^{\infty}} \Bigg \}, \quad
m+1\leq k\leq q.
\end{equation}
\begin{equation}
\label{Pk3}
P_{k}:=\hbox{max}\Bigg\{ \Bigg\|\frac{G_{k, \ n}}{|n|}\Bigg\|_{l^{\infty}},
\Bigg\|nG_{k, \ n}\Bigg\|_{l^{\infty}} \Bigg \}, \quad q+1\leq k\leq N_{1}.
\end{equation}
\begin{equation}
\label{Pk4}
P_{k}:=\hbox{max}\Bigg\{ \Bigg\|\frac{G_{k, \ n}}{|n|+a_{k}}\Bigg\|_{l^{\infty}},
\Bigg\|\frac{n^{2}G_{k, \ n}}{|n|+a_{k}}\Bigg\|_{l^{\infty}} \Bigg \}, \quad
N_{1}+1\leq k\leq N_{2}.
\end{equation}
We recall Lemma A3 of \cite{VV18}. Hence, under Assumption 1.4 expressions
(\ref{Pk1}), (\ref{Pk2}) and (\ref{Pk3}) are finite.

\noindent
It can be easily verified that (\ref{Pk4}) are finite as well. Clearly,
for $N_{1}+1\leq k \leq N_{2}$ by means of (\ref{fiub}) we derive
$$
\Bigg|\frac{G_{k, \ n}}{|n|+a_{k}}\Bigg|\leq \frac{|G_{k, \ n}|}{a_{k}}\leq
\frac{1}{a_{k}\sqrt{2\pi}}\|G_{k}\|_{L^{1}(I)}<\infty
$$
due the one of our assumptions. Let us use (\ref{fiub}) to obtain
$$
\Bigg|\frac{n^{2}G_{k, \ n}}{|n|+a_{k}}\Bigg|\leq |nG_{k, \ n}|\leq
\frac{1}{\sqrt{2\pi}}\|G_{k}'(x)\|_{L^{1}(I)}<\infty.
$$
Therefore, $P_{k}<\infty$ for $N_{1}+1\leq k\leq N_{2}$ as well. This enables us
to define
\begin{equation}
\label{p}  
P:=\hbox{max}P_{k}, \quad 1\leq k\leq N_{2}
\end{equation}
with $P_{k}$ given by formulas (\ref{Pk1}), (\ref{Pk2}), (\ref{Pk3})
and (\ref{Pk4}).

\noindent
To study the existence of stationary solutions of our
system of equations, we use 
the corresponding function space 
$$
H^{2}(I):=\{v(x):I\to {\mathbb R} \ | \ v(x), v''(x)\in L^{2}(I), \quad
v(0)=v(2\pi), \quad v'(0)=v'(2\pi) \}.
$$
Let us aim at $u_{k}(x)\in H^{2}(I), \ 1\leq k\leq m$ and
$N_{1}+1\leq k\leq N_{2}$ as well.

\noindent
We introduce the following auxiliary constrained subspaces
$$
H_{k}^{2}(I):=\Bigg\{v\in H^{2}(I) \ | \ \Big(v(x), \frac{e^{\pm in_{k}x}}
{\sqrt{2\pi}}\Big)_{L^{2}(I)}=0 \Bigg\}, \quad n_{k}\in {\mathbb N}, \quad
m+1\leq k\leq q.
$$
Our goal is to have $u_{k}(x)\in H_{k}^{2}(I), \ m+1\leq k\leq q$. Similarly,
$$
H_{0}^{2}(I):=\{v\in H^{2}(I) \ | \ (v(x),1)_{L^{2}(I)}=0 \}.
$$
The goal is to have $u_{k}(x)\in H_{0}^{2}(I), \ q+1\leq k\leq N_{1}$. The 
constrained subspaces mentioned above are Hilbert spaces as well 
(see e.g. Chapter 2.1 of ~\cite{HS96}).

\noindent
The resulting space used to
establish the existence of solutions $u(x): I\to {\mathbb R}^{N_{2}}$ 
of system (\ref{ae1}), (\ref{ae2}) will be the direct sum of the spaces
introduced above, namely
$$
H_{c}^{2}(I, \ {\mathbb R}^{N_{2}}):=\oplus_{k=1}^{m}H^{2}(I)\oplus_{k=m+1}^{q}
H_{k}^{2}(I)\oplus_{k=q+1}^{N_{1}}H_{0}^{2}(I)\oplus_{k=N_{1}+1}^{N_{2}}H^{2}(I).
$$
The corresponding Sobolev norm equals to
$$
\|u\|_{H_{c}^{2}(I, \ {\mathbb R}^{N_{2}})}^{2}:=\sum_{k=1}^{N_{2}}\{\|u_{k}\|_{L^{2}(I)}^{2}+
\|u_{k}''\|_{L^{2}(I)}^{2}\},
$$
where $u(x):I \to {\mathbb R}^{N_{2}}$.

\noindent
Let us demonstrate that the system of equations (\ref{ae1}), (\ref{ae2}) in
this case defines a map on the space stated above, which will be a strict
contraction under the given conditions.

\bigskip

\noindent
{\bf Theorem 1.5.} {\it Let $\Omega=I$, Assumptions 1.1 and 1.4 are valid and
$2\sqrt{\pi}PL<1$.   

\medskip
  
\noindent  
Then the map $\tau_{a}v=u$ on $H_{c}^{2}(I, \ {\mathbb R}^{N_{2}})$
defined by the system of equations (\ref{ae1}), (\ref{ae2}) has a unique fixed
point $v_{a}(x): I\to {\mathbb R}^{N_{2}}$. This is the
only stationary solution of system (\ref{id1}), (\ref{id11}) in 
$H_{c}^{2}(I, \ {\mathbb R}^{N_{2}})$.

\medskip

\noindent  
This fixed point $v_{a}(x)$ does not vanish identically in $I$
provided the Fourier coefficients 
$G_{k, \ n}F_{k}(0, x)_{n}\neq 0$ for some $1\leq k\leq N_{2}$ and a certain
$n\in {\mathbb Z}$.}

\bigskip

\noindent
Let us note that the the Fredholm operators
$$
\sqrt{-{d^{2}\over dx^{2}}}-n_{k}: H_{k}^{1}(I)\to L^{2}(I) \quad and \quad
\sqrt{-{d^{2}\over dx^{2}}}: H_{0}^{1}(I)\to L^{2}(I)
$$
involved in our system (\ref{ae1}), (\ref{ae2}) with $\Omega=I$
have the trivial kernels. Here
$$
H_{k}^{1}(I):=\Bigg\{v\in H^{1}(I) \ | \ \Big(v(x), \frac{e^{\pm in_{k}x}}
{\sqrt{2\pi}}\Big)_{L^{2}(I)}=0 \Bigg\}, \quad n_{k}\in {\mathbb N}, \quad
m+1\leq k\leq q
$$
and
$$
H_{0}^{1}(I):=\{v\in H^{1}(I) \ | \ (v(x),1)_{L^{2}(I)}=0 \}.
$$

\bigskip

\noindent
Let us conclude the work with the studies of our system of equations in 
the layer domain. This is the product of the two sets, so that one is 
the $I$ interval with periodic boundary conditions as in the previous part 
of the article and another is the whole space of dimension either one or two,
namely
$\Omega=I\times {\mathbb R}^{d}=[0, 2\pi]\times {\mathbb R}^{d}, \ d=1,2$ and
$x=(x_{1}, x_{\perp})$, where $x_{1}\in I$ and
$x_{\perp}\in {\mathbb R}^{d}$. The resulting Laplacian in this context
will be $\displaystyle{\Delta:=\frac{\partial^{2}}{\partial x_{1}^{2}}+
\sum_{s=1}^{d}\frac{\partial^{2}}{\partial x_{\perp, \ s}^{2}}}$.

\noindent
The corresponding Sobolev space for this problem will be 
$H^{2}(\Omega, \ {\mathbb R}^{N_{2}})$  of vector
functions $u(x): \Omega\to {\mathbb R}^{N_{2}}$, such that for $1\leq k\leq N_{2}$
$$
u_{k}(x), \ \Delta u_{k}(x)\in 
L^{2}(\Omega), \quad u_{k}(0,x_{\perp})=u_{k}(2\pi,x_{\perp}), \quad \frac
{\partial u_{k}}{\partial x_{1}}(0, x_{\perp})=\frac{\partial u_{k}}
{\partial x_{1}}(2\pi, x_{\perp}),
$$
where $x_{\perp}\in {\mathbb R}^{d}$. It is equipped with the norm
$$
\|u\|_ {H^{2}(\Omega, \ {\mathbb R}^{N_{2}})}^{2}=\sum_{k=1}^{N_{2}}
\{\|u_{k}\|_ {L^{2}(\Omega)}^{2}+\|\Delta u_{k}\|_ {L^{2}(\Omega)}^{2}\}.
$$
Similarly to
the whole space case treated in Theorem 1.3, the operators
$\sqrt{-\Delta}-a_{k}: H^{1}(\Omega)\to L^{2}(\Omega)$ for
$a_{k}\geq 0, \ 1\leq k\leq N_{1}$ do not possess 
the Fredholm property.

\noindent
We establish that problem (\ref{ae1}), (\ref{ae2})
in this case defines a map
$t_{a}:H^{2}(\Omega, \ {\mathbb R}^{N_{2}})\to H^{2}(\Omega, \ {\mathbb R}^{N_{2}})$, 
which is a strict contraction under the appropriate technical assumptions
formulated below.

\bigskip

\noindent
{\bf Assumption 1.6.} {\it Let $G_{k}(x): \Omega \to {\mathbb R}, \ G_{k}(x)\in
C(\Omega)\cap W^{1, 1}(\Omega), \  G_{k}(0, x_{\perp})=G_{k}(2\pi, x_{\perp})$ and 
$F_{k}(u,0, x_{\perp})=F_{k}(u,2\pi, x_{\perp})$ for $x_{\perp}\in {\mathbb R}^{d}$, 
$u\in {\mathbb R}^{N_{2}}, \ d=1,2$ and $1\leq k\leq N_{2}$. Let
$N_{1}\geq 3, \ 1\leq m<q\leq N_{1}-1$ with $m,q \in {\mathbb N}$ and
$N_{2}>N_{1}$.

\medskip

\noindent
I) Assume for  $1\leq k\leq m$ we have $n_{k}<a_{k}<n_{k}+1, \ n_{k}\in 
{\mathbb Z}^{+}={\mathbb N}\cup \{0\}, \ x_{\perp}G_{k}(x)\in L^{1}(\Omega)$ and
\begin{equation}
\label{or6}
\Bigg(G_{k}(x_{1},x_{\perp}),\frac{e^{inx_{1}}}{\sqrt{2\pi}}\frac
{e^{\pm i\sqrt{a_{k}^{2}-n^{2}}x_{\perp}}}{\sqrt{2\pi}}\Bigg)_{L^{2}(\Omega)}=0, \quad |n|\leq 
n_{k} \quad when \quad d=1,
\end{equation}
\begin{equation}
\label{or7}
\Bigg(G_{k}(x_{1},x_{\perp}),\frac{e^{inx_{1}}}{\sqrt{2\pi}}\frac
{e^{ipx_{\perp}}}{2\pi}\Bigg)_{L^{2}(\Omega)}=0, \quad p\in S_{\sqrt{a_{k}^{2}-n^{2}}}^{2}
\quad, \quad |n|\leq n_{k} \quad for \quad d=2.
\end{equation}

\medskip

\noindent
II) Assume for $m+1\leq k\leq q$ we have 
$a_{k}=n_{k}, \ n_{k}\in {\mathbb N}, \ x_{\perp}^{2}G_{k}(x)\in L^{1}(\Omega)$
and
\begin{equation}
\label{or8}
\Bigg(G_{k}(x_{1},x_{\perp}),\frac{e^{inx_{1}}}{\sqrt{2\pi}}\frac
{e^{\pm i\sqrt{n_{k}^{2}-n^{2}}x_{\perp}}}{\sqrt{2\pi}}\Bigg)_{L^{2}(\Omega)}=0, \quad |n|\leq 
n_{k}-1 \ when \ d=1,
\end{equation}
\begin{equation}
\label{or9}
\Bigg(G_{k}(x_{1},x_{\perp}),\frac{e^{inx_{1}}}{\sqrt{2\pi}}\frac
{e^{ipx_{\perp}}}{2\pi}\Bigg)_{L^{2}(\Omega)}=0, \quad p\in S_{\sqrt{n_{k}^{2}-n^{2}}}^{2},
\quad  |n|\leq n_{k}-1 \quad for \quad d=2,
\end{equation}
\begin{equation}
\label{or10}
\Bigg(G_{k}(x_{1},x_{\perp}),\frac{e^{\pm in_{k}x_{1}}}{\sqrt{2\pi}}\Bigg)_{L^{2}(\Omega)}
=0, \quad \Bigg(G_{k}(x_{1},x_{\perp}),\frac{e^{\pm in_{k}x_{1}}}{\sqrt{2\pi}}
x_{\perp, \ s}\Bigg)_{L^{2}(\Omega)}=0, 
\end{equation}
where $1\leq s\leq d$.

\medskip

\noindent
III) Assume for $q+1\leq k\leq N_{1}$ we have $a_{k}=0, \quad x_{\perp}G_{k}(x)\in
L^{1}(\Omega)$ and
\begin{equation}
\label{or11}
(G_{k}(x),1)_{L^{2}(\Omega)}=0.
\end{equation}

\medskip

\noindent
IV) Let $a_{k}>0, \ N_{1}+1\leq k\leq N_{2}$.
}

\bigskip

\noindent
Let us make use of the Fourier transform for the functions on such a product
of sets, namely
\begin{equation}
\label{ftl}  
{\widehat{G_{k, \ n}}}(p):=\frac{1}{(2\pi)^{\frac{d+1}{2}}}\int_{{\mathbb R}^{d}}dx_{\perp}
e^{-ipx_{\perp}}\int_{0}^{2\pi}G_{k}(x_{1},x_{\perp})e^{-inx_{1}}dx_{1},
\end{equation}
where  $p\in {\mathbb R}^{d}, \ n\in {\mathbb Z}, \ 1\leq k\leq N_{2}$.
Obviously,
\begin{equation}
\label{fubl}  
\|\widehat{G_{k, \ n}}(p)\|_{L^{\infty}_{n, p}}:=\hbox{sup}_{\{p\in {\mathbb R}^{d},
\quad n\in {\mathbb Z}\}}|\widehat{G_{k, \ n}}(p)|\leq
\frac{1}{(2\pi)^{\frac{d+1}{2}}}\|G_{k}\|_{L^{1}(\Omega)},
\end{equation}
\begin{equation}
\label{fubld}
\|n\widehat{G_{k, \ n}}(p)\|_{L^{\infty}_{n, p}}\leq \frac{1}{(2\pi)^{\frac{d+1}{2}}}
\Big\|\frac{\partial G_{k}}{\partial x_{1}}\Big\|_{L^{1}(\Omega)}, \quad
\|p\widehat{G_{k, \ n}}(p)\|_{L^{\infty}_{n, p}}\leq \frac{1}{(2\pi)^{\frac{d+1}{2}}}
\|\nabla_{x_{\perp}}G_{k}\|_{L^{1}(\Omega)}.
\end{equation}
Here $\nabla_{x_{\perp}}$ stands for the gradient with respect to $x_{\perp}$.

\noindent
We define the following technical expressions
\begin{equation}
\label{R1}
R_{k}:=\hbox{max}\Bigg\{\Bigg\|\frac{\widehat{G_{k, \ n}}(p)}
{\sqrt{p^{2}+n^{2}}-a_{k}}
\Bigg\|_{L^{\infty}_{n, p}}, \ \Bigg\|\frac{(p^{2}+n^{2})\widehat{G_{k, \ n}}(p)}
{\sqrt{p^{2}+n^{2}}-a_{k}}\Bigg\|_{L^{\infty}_{n, p}}\Bigg\}, \quad 1\leq k\leq m.
\end{equation}
\begin{equation}
\label{R2}
R_{k}:=\hbox{max}\Bigg\{\Bigg\|\frac{\widehat{G_{k, \ n}}(p)}
{\sqrt{p^{2}+n^{2}}-n_{k}}
\Bigg\|_{L^{\infty}_{n, p}}, \ \Bigg\|\frac{(p^{2}+n^{2})\widehat{G_{k, \ n}}(p)}
{\sqrt{p^{2}+n^{2}}-n_{k}}\Bigg\|_{L^{\infty}_{n, p}}\Bigg\}, \quad m+1\leq k\leq q.
\end{equation}
\begin{equation}
\label{R3}
R_{k}:=\hbox{max}\Bigg\{\Bigg\|\frac{\widehat{G_{k, \ n}}(p)}{\sqrt{p^{2}+n^{2}}}
\Bigg\|_{L^{\infty}_{n, p}}, \ \Bigg\|\sqrt{p^{2}+n^{2}}\widehat{G_{k, \ n}}(p)
\Bigg\|_{L^{\infty}_{n, p}}\Bigg\}, \quad q+1\leq k\leq N_{1}.
\end{equation}
\begin{equation}
\label{R4}
R_{k}:=\hbox{max}\Bigg\{\Bigg\|\frac{\widehat{G_{k, \ n}}(p)}
{\sqrt{p^{2}+n^{2}}+a_{k}}
\Bigg\|_{L^{\infty}_{n, p}}, \ \Bigg\|\frac{(p^{2}+n^{2})\widehat{G_{k, \ n}}(p)}
{\sqrt{p^{2}+n^{2}}+a_{k}}\Bigg\|_{L^{\infty}_{n, p}}\Bigg\}, \quad
N_{1}+1\leq k\leq N_{2}.
\end{equation}
Let us recall Lemmas A4, A5 and A6 of ~\cite{VV18}.
Thus, Assumption 1.6 implies that
the quantities given by (\ref{R1}), (\ref{R2}) and (\ref{R3}) are
finite.

\noindent
It can be trivially checked that (\ref{R4}) are finite as well.
Evidently, for $N_{1}+1\leq k\leq N_{2}$ by virtue of (\ref{fubl}) we have
$$
\Big|\frac{\widehat{G_{k, \ n}}(p)}{\sqrt{p^{2}+n^{2}}+a_{k}}\Big|\leq
\frac{|{\widehat{G_{k, \ n}}(p)}|}{a_{k}}\leq \frac{1}{(2\pi)^{\frac{d+1}{2}}a_{k}}
\|G_{k}\|_{L^{1}(\Omega)}<\infty
$$
via the one of our assumptions. We use (\ref{fubld}) to derive that
$$
\Big|\frac{(p^{2}+n^{2})\widehat{G_{k, \ n}}(p)}{\sqrt{p^{2}+n^{2}}+a_{k}}\Big|\leq
\sqrt{p^{2}+n^{2}}|\widehat{G_{k, \ n}}(p)|\leq
\|p\widehat{G_{k, \ n}}(p)\|_{L^{\infty}_{n, p}}+\|n\widehat{G_{k, \ n}}(p)\|_{L^{\infty}_{n, p}}
\leq
$$
$$
\frac{1}{(2\pi)^{\frac{d+1}{2}}}
\|\nabla_{x_{\perp}}G_{k}\|_{L^{1}(\Omega)}+\frac{1}{(2\pi)^{\frac{d+1}{2}}}
\Big\|\frac{\partial G_{k}}{\partial x_{1}}\Big\|_{L^{1}(\Omega)}<\infty
$$
as well. This means that $R_{k}<\infty$ when $N_{1}+1\leq k\leq N_{2}$.
This allowes us to introduce
\begin{equation}
\label{R}  
R:=\hbox{max}R_{k}, \quad 1\leq k\leq N_{2},
\end{equation}
where $R_{k}$ are defined in (\ref{R1}), (\ref{R2}), (\ref{R3}) and (\ref{R4}).
Our final statement is as follows.

\bigskip

\noindent
{\bf Theorem 1.7.} {\it Let $\Omega=I\times {\mathbb R}^{d}, \ d=1,2$, 
Assumptions 1.1 and 1.6 hold and $\sqrt{2}(2\pi)^{\frac{d+1}{2}}RL<1$.

\medskip

\noindent  
Then the map $t_{a}v=u$ on $H^{2}(\Omega, \ {\mathbb R}^{N_{2}})$ defined 
by system (\ref{ae1}), (\ref{ae2}) has a unique fixed point
$v_{a}(x):\Omega\to {\mathbb R}^{N_{2}}$. This is the only stationary solution
of the system of equations (\ref{id1}), (\ref{id11}) in
$H^{2}(\Omega, \ {\mathbb R}^{N_{2}})$.

\medskip

\noindent
This fixed point $v_{a}(x)$ is nontrivial provided that the intersection of
supports of the Fourier images of the functions
$supp\widehat{F_{k}(0,x)}_{n}(p)\cap supp \widehat{G_{k, \ n}}(p)$
is a set of nonzero Lebesgue measure in ${\mathbb R}^{d}$ for a certain
$1\leq k\leq N_{2}$ and some $n\in {\mathbb Z}$.}

\bigskip

\noindent
Note that the maps used in the theorems above are applied to the
real valued vector functions by means of the conditions on $F_{k}(u,x)$ and 
$G_{k}(x), \ 1\leq k\leq N_{2}$ contained in the nonlocal terms of
(\ref{ae1}), (\ref{ae2}).


\setcounter{equation}{0}

\section{The System in the Whole Space}

\bigskip

{\it Proof of Theorem 1.3.} First we suppose that in the situation when
$\Omega={\mathbb R}^{d}, \ 1\leq d\leq 3$ there exists 
$v(x)\in H^{2}({\mathbb R}^{d}, \ {\mathbb R}^{N_{2}})$, so that the system of
equations (\ref{ae1}), (\ref{ae2}) has two solutions 
$u^{(1),(2)}(x)\in H^{2}({\mathbb R}^{d}, \ {\mathbb R}^{N_{2}})$.
Then the difference vector function 
$w(x):=u^{(1)}(x)-u^{(2)}(x)\in H^{2}({\mathbb R}^{d}, \ {\mathbb R}^{N_{2}})$ 
solves the homogeneous system 
$$
\sqrt{-\Delta} w_{k}=a_{k}w_{k}, \quad a_{k}\geq 0, \quad 1\leq k\leq N_{1},
$$
$$
\sqrt{-\Delta} w_{k}=-a_{k}w_{k}, \quad a_{k}>0, \quad N_{1}+1\leq k\leq N_{2}.
$$
Because the operator
$\sqrt{-\Delta}: H^{1}({\mathbb R}^{d})\to L^{2}({\mathbb R}^{d})$
does not possess any nontrivial eigenfunctions, we obtain that $w_{k}(x)$
vanish in ${\mathbb R}^{d}$ for all $1\leq k\leq N_{2}$. 

\noindent
Let us consider an arbitrary vector function 
$v(x)\in H^{2}({\mathbb R}^{d}, \ {\mathbb R}^{N_{2}})$  and apply the standard
Fourier transform (\ref{ft}) to both sides of the system of equations
(\ref{ae1}), (\ref{ae2}). This yields
\begin{equation}
\label{f1}
\widehat{u_{k}}(p)=(2\pi)^{\frac{d}{2}}\frac{\widehat{G_{k}}(p)\widehat{f_{k}}(p)}
{|p|-a_{k}}, \quad a_{k}\geq 0, \quad 1\leq k\leq N_{1},
\end{equation}
\begin{equation}
\label{f11}
\widehat{u_{k}}(p)=(2\pi)^{\frac{d}{2}}\frac{\widehat{G_{k}}(p)\widehat{f_{k}}(p)}
{|p|+a_{k}}, \quad a_{k}>0, \quad N_{1}+1\leq k\leq N_{2}.
\end{equation}
Here $\widehat{f_{k}}(p)$ stands for the Fourier transform of $F_{k}(v(x),x)$.
Also,
\begin{equation}
\label{f12}
p^{2}\widehat{u_{k}}(p)=(2\pi)^{\frac{d}{2}}\frac{p^{2}\widehat{G_{k}}(p)
\widehat{f_{k}}(p)}{|p|-a_{k}}, \quad a_{k}\geq 0, \quad 1\leq k\leq N_{1},
\end{equation}
\begin{equation}
\label{f112}
p^{2}\widehat{u_{k}}(p)=(2\pi)^{\frac{d}{2}}\frac{p^{2}\widehat{G_{k}}(p)
\widehat{f_{k}}(p)}{|p|+a_{k}}, \quad a_{k}>0, \quad N_{1}+1\leq k\leq N_{2}.
\end{equation}

\noindent
We derive the following upper bounds using quantities (\ref{Mka}), (\ref{Mk0})
(\ref{Mka+}) and (\ref{M}), namely
$$
|\widehat{u_{k}}(p)|\leq (2\pi)^{d\over 2}M|\widehat{f_{k}}(p)|  \quad and
\quad |p^{2}\widehat{u_{k}}(p)|\leq (2\pi)^{d\over 2}M|\widehat{f_{k}}(p)|, \
1\leq k \leq N_{2}.
$$
This enables us to obtain the estimate from above for the norm as
$$  
\|u\|_{H^{2}({\mathbb R}^{d}, \ {\mathbb R}^{N_{2}})}^{2}=\sum_{k=1}^{N_{2}}
\{\|\widehat{u_{k}}(p)\|_{L^{2}({\mathbb R}^{d})}^{2}+
\|p^{2}\widehat{u_{k}}(p)\|_{L^{2}({\mathbb R}^{d})}^{2}\}=
$$
$$
\sum_{k=1}^{N_{2}}\Big\{\int_{{\mathbb R}^{d}}|\widehat{u_{k}}(p)|^{2}dp+
\int_{{\mathbb R}^{d}}|p^{2}\widehat{u_{k}}(p)|^{2}dp\Big\}\leq
$$
\begin{equation}
\label{uh2}
2(2\pi)^{d}M^{2}\sum_{k=1}^{N_{2}}\int_{{\mathbb R}^{d}}|\widehat{f_{k}}(p)|^{2}dp=  
2(2\pi)^{d}M^{2}\sum_{k=1}^{N_{2}}\|F_{k}(v(x),x)\|_{L^{2}({\mathbb R}^{d})}^{2}.
\end{equation}
Let us recall inequality (\ref{ub1}) of Assumption 1.1 above. Thus,
all $F_{k}(v(x),x)\in L^{2}({\mathbb R}^{d})$, such that 
the right side of (\ref{uh2}) is finite.

\noindent
Hence, for any
$v(x)\in H^{2}({\mathbb R}^{d}, \ {\mathbb R}^{N_{2}})$ there exists a unique 
vector function $u(x)\in H^{2}({\mathbb R}^{d}, \ {\mathbb R}^{N_{2}})$
satistying system (\ref{ae1}), (\ref{ae2}). Its Fourier image is given by
(\ref{f1}), (\ref{f11}). Hence, the map
$T_{a}:  H^{2}({\mathbb R}^{d}, \ {\mathbb R}^{N_{2}})\to  
H^{2}({\mathbb R}^{d}, \ {\mathbb R}^{N_{2}})$ is well defined.

\noindent
This allows us to choose arbitrarily
$v^{(1),(2)}(x)\in H^{2}({\mathbb R}^{d}, \ {\mathbb R}^{N_{2}})$ and obtain
their images under such map 
$u^{(1),(2)}:=T_{a}v^{(1),(2)}\in H^{2}({\mathbb R}^{d}, \ {\mathbb R}^{N_{2}})$.
By virtue of (\ref{ae1}), (\ref{ae2}),
\begin{equation}
\label{ae11}
\sqrt{-\Delta} u_{k}^{(1)}-a_{k}u_{k}^{(1)}=\int_{{\mathbb R}^{d}}G_{k}(x-y)F_{k}
(v_{1}^{(1)}(y),v_{2}^{(1)}(y),..., v_{N_{2}}^{(1)}(y), y)dy, \quad a_{k}\geq 0
\end{equation}
if $1\leq k\leq N_{1}$,
\begin{equation}
\label{ae21}
\sqrt{-\Delta} u_{k}^{(1)}+a_{k}u_{k}^{(1)}=\int_{{\mathbb R}^{d}}G_{k}(x-y)
F_{k}(v_{1}^{(1)}(y),v_{2}^{(1)}(y),..., v_{N_{2}}^{(1)}(y), y)dy, \quad a_{k}>0
\end{equation}
if $N_{1}+1\leq k\leq N_{2}$. Analogously,
\begin{equation}
\label{ae12}
\sqrt{-\Delta} u_{k}^{(2)}-a_{k}u_{k}^{(2)}=\int_{{\mathbb R}^{d}}G_{k}(x-y)F_{k}
(v_{1}^{(2)}(y),v_{2}^{(2)}(y),..., v_{N_{2}}^{(2)}(y), y)dy, \quad a_{k}\geq 0
\end{equation}
for $1\leq k\leq N_{1}$,
\begin{equation}
\label{ae22}
\sqrt{-\Delta} u_{k}^{(2)}+a_{k}u_{k}^{(2)}=\int_{{\mathbb R}^{d}}G_{k}(x-y)
F_{k}(v_{1}^{(2)}(y),v_{2}^{(2)}(y),..., v_{N_{2}}^{(2)}(y), y)dy, \quad a_{k}>0
\end{equation}
for $N_{1}+1\leq k\leq N_{2}$.

\noindent
We apply the standard Fourier transform
(\ref{ft}) to both sides of the systems of equations
(\ref{ae11}), (\ref{ae21}) and (\ref{ae12}), (\ref{ae22}). This gives us
\begin{equation}
\label{uk1hp}
\widehat{u_{k}^{(1)}}(p)=(2\pi)^{\frac{d}{2}}
\frac{\widehat{G_{k}}(p)\widehat{f_{k}^{(1)}}(p)}{|p|-a_{k}}, \quad a_{k}\geq 0,
\quad 1\leq k\leq N_{1},
\end{equation}  
\begin{equation}
\label{uk1hp2}
\widehat{u_{k}^{(1)}}(p)=(2\pi)^{\frac{d}{2}}
\frac{\widehat{G_{k}}(p)\widehat{f_{k}^{(1)}}(p)}{|p|+a_{k}}, \quad a_{k}>0, \quad
N_{1}+1\leq k\leq N_{2}
\end{equation} 
and
\begin{equation}
\label{uk2hp}
\widehat{u_{k}^{(2)}}(p)=(2\pi)^{\frac{d}{2}}
\frac{\widehat{G_{k}}(p)\widehat{f_{k}^{(2)}}(p)}{|p|-a_{k}}, \quad a_{k}\geq 0,
\quad 1\leq k\leq N_{1},
\end{equation}  
\begin{equation}
\label{uk2hp2}
\widehat{u_{k}^{(2)}}(p)=(2\pi)^{\frac{d}{2}}
\frac{\widehat{G_{k}}(p)\widehat{f_{k}^{(2)}}(p)}{|p|+a_{k}}, \quad a_{k}>0, \quad
N_{1}+1\leq k\leq N_{2}.
\end{equation} 
Here $\widehat{f_{k}^{(1),(2)}}(p)$ denote the Fourier transforms of
$F_{k}(v^{(1),(2)}(x),x)$.

\noindent
Using formulas (\ref{uk1hp}), (\ref{uk1hp2}),
(\ref{uk2hp}), (\ref{uk2hp2}) along with (\ref{Mka}), (\ref{Mk0})
(\ref{Mka+}) and (\ref{M}), we easily derive the estimates from above for
$1\leq k \leq N_{2}$ as
$$
\Bigg|\widehat{u_{k}^{(1)}}(p)-\widehat{u_{k}^{(2)}}(p)\Bigg|\leq (2\pi)^{d\over 2}M
\Bigg|\widehat{f_{k}^{(1)}}(p)-\widehat{f_{k}^{(2)}}(p)\Bigg|,
$$
$$
\Bigg|p^{2}\widehat{u_{k}^{(1)}}(p)-p^{2}\widehat{u_{k}^{(2)}}(p)\Bigg|\leq 
(2\pi)^{d\over 2}M\Bigg|\widehat{f_{k}^{(1)}}(p)-\widehat{f_{k}^{(2)}}(p)\Bigg|.
$$
This enables us to obtain the upper bound on the corresponding 
norm of the difference of vector functions 
$$
\|u^{(1)}-u^{(2)}\|_{H^{2}({\mathbb R}^{d}, \ {\mathbb R}^{N_{2}})}^{2}=
\sum_{k=1}^{N_{2}}\Big\{\Big\|\widehat{u_{k}^{(1)}}(p)-\widehat{u_{k}^{(2)}}(p)\Big\|
_{L^{2}({\mathbb R}^{d})}^{2}+\Big
\|p^{2}\Big[\widehat{u_{k}^{(1)}}(p)-\widehat{u_{k}^{(2)}}(p)\Big]\Big\|
_{L^{2}({\mathbb R}^{d})}^{2}\Big\}=
$$
$$
\sum_{k=1}^{N_{2}}\Big\{\int_{{\mathbb R}^{d}}\Big|\widehat{u_{k}^{(1)}}(p)-
\widehat{u_{k}^{(2)}}(p)\Big|^{2}dp+\int_{{\mathbb R}^{d}}
\Big|p^{2}\widehat{u_{k}^{(1)}}(p)-p^{2}
\widehat{u_{k}^{(2)}}(p)\Big|^{2}dp\Big\}\leq
$$
$$
2(2\pi)^{d}M^{2}\sum_{k=1}^{N_{2}}\int_{{\mathbb R}^{d}}\Big|\widehat{f_{k}^{(1)}}(p)-
\widehat{f_{k}^{(2)}}(p)\Big|^{2}dp=
$$
\begin{equation}
\label{u12np}
2(2\pi)^{d}M^{2}
\sum_{k=1}^{N_{2}}\|F_{k}(v^{(1)}(x),x)-F_{k}(v^{(2)}(x),x)\|_{L^{2}({\mathbb R}^{d})}^{2}.
\end{equation}
Note that by virtue of the Sobolev embedding theorem, for $1\leq k\leq N_{2}$
we have 
$v^{(1),(2)}_{k}(x)\in H^{2}({\mathbb R}^{d})\subset
L^{\infty}({\mathbb R}^{d}), \ 1\leq d\leq 3$.

\noindent
We recall inequality (\ref{lk1}) of Assumption 1.1 above. This yields
$$
\sum_{k=1}^{N_{2}}\|F_{k}(v^{(1)}(x),x)-F_{k}(v^{(2)}(x),x)\|_{L^{2}({\mathbb R}^{d})}^{2}\leq
L^{2}\sum_{k=1}^{N_{2}}\|v_{k}^{(1)}(x)-v_{k}^{(2)}(x)\|_{L^{2}({\mathbb R}^{d})}^{2},
$$
such that
\begin{equation}
\label{fp1}  
\|T_{a}v^{(1)}-T_{a}v^{(2)}\|_{H^{2}({\mathbb R}^{d}, \ {\mathbb R}^{N_{2}})}\leq
\sqrt{2}(2\pi)^{\frac{d}{2}}
ML\|v^{(1)}-v^{(2)}\|_{H^{2}({\mathbb R}^{d}, \  {\mathbb R}^{N_{2}})}.
\end{equation}
The constant in the right side of bound (\ref{fp1}) is less than one we assume. 
Hence, the Fixed Point Theorem yields the 
existence of a unique vector function
$v_{a}(x)\in H^{2}({\mathbb R}^{d}, {\mathbb R}^{N_{2}})$, such that
$T_{a}v_{a}=v_{a}$.
This is the only stationary solution of the system of equations
(\ref{id1}), (\ref{id11}) in
$H^{2}({\mathbb R}^{d}, \ {\mathbb R}^{N_{2}})$.

\noindent
Let us suppose that
$v_{a}(x)$ is trivial in the whole space. 
This will give us the contradiction to the condition that for some
$1\leq k \leq N_{2}$ the Fourier images of $G_{k}(x)$ and $F_{k}(0,x)$ do
not vanish simultaneously on a certain set of nonzero Lebesgue measure in
${\mathbb R}^{d}$.  \hfill\lanbox


 \setcounter{equation}{0}

 \section{The System on the $[0, 2\pi]$ Interval}

\bigskip

{\it Proof of Theorem 1.5.}  First we suppose that for a certain
$v(x)\in H_{c}^{2}(I, \ {\mathbb R}^{N_{2}})$ there exist two
solutions $u^{(1),(2)}(x)\in H_{c}^{2}(I, \ {\mathbb R}^{N_{2}})$ 
of the system of equations (\ref{ae1}), (\ref{ae2}) with $\Omega=I$.
Then the difference vector function 
$w(x):=u^{(1)}(x)-u^{(2)}(x)\in  H_{c}^{2}(I, \ {\mathbb R}^{N_{2}})$ will satisfy
the homogeneous system 
\begin{equation}
\label{his1}  
\sqrt{-\frac{d^{2}}{dx^{2}}}w_{k}=a_{k}w_{k}, \quad a_{k}\geq 0, \quad 1\leq k
\leq N_{1},
\end{equation}
\begin{equation}
\label{his2} 
\sqrt{-\frac{d^{2}}{dx^{2}}}w_{k}=-a_{k}w_{k}, \quad a_{k}>0, \quad N_{1}+1\leq k
\leq N_{2}.
\end{equation}
Let us recall Assumption 1.4. Hence,
$a_{k}>0, \ a_{k}\neq n, \ n\in {\mathbb N}$ for $1\leq k\leq m$. Thus,  
they are not the eigenvalues of the operator
\begin{equation}
\label{iop}  
\sqrt{-{d^{2}\over dx^{2}}}: H^{1}(I)\to L^{2}(I).
\end{equation}
Therefore, $w_{k}(x)$ vanish identically in $I$ for $1\leq k\leq m$.

\noindent
For $m+1\leq k \leq q$ the values of $a_{k}$ coincide with the eigenvalues of
(\ref{iop}) as we assume.
But $w_{k}$ are contained in the constraned subspaces
$H_{k}^{2}(I)$. Hence, $w_{k}$ are trivial in $I$ when
$m+1\leq k \leq q$  due to their
orthogonality to the eigenfunctions 
$\displaystyle{\frac{e^{\pm in_{k}x}}{\sqrt{2\pi}}}$ of (\ref{iop}).

\noindent
By means of Assumption 1.4, 
the constants $a_{k}$ vanish for $q+1\leq k\leq N_{1}$. But
$w_{k}\in H_{0}^{2}(I)$. Thus, they are orthogonal to the zero mode of
operator (\ref{iop}). This means that $w_{k}(x)$ are trivial in $I$ as well for
$q+1\leq k\leq N_{1}$.

\noindent
Finally, we deal with the case when
$N_{1}+1\leq k\leq N_{2}$. Clearly, operator (\ref{iop}) cannot have any negative
eigenvalues. By virtue of (\ref{his2}) we obtain that $w_{k}(x)$
vanish identically in $I$ for $N_{1}+1\leq k\leq N_{2}$. 

\noindent
Let us choose arbitrarily $v(x)\in H_{c}^{2}(I, \ {\mathbb R}^{N_{2}})$. We 
apply the Fourier transform (\ref{fti}) to both sides of system 
(\ref{ae1}),  (\ref{ae2}) studied on the interval $[0, 2\pi]$.
This gives us
\begin{equation}
\label{f2}
u_{k, \ n}=\sqrt{2\pi}\frac{G_{k, \ n}f_{k, \ n}}{|n|-a_{k}}, \quad a_{k}\geq 0, 
\quad 1\leq k\leq N_{1}, \quad n\in {\mathbb Z},
\end{equation}
\begin{equation}
\label{f21}
u_{k, \ n}=\sqrt{2\pi}\frac{G_{k, \ n}f_{k, \ n}}{|n|+a_{k}}, \quad a_{k}>0, 
\quad N_{1}+1\leq k\leq N_{2}, \quad n\in {\mathbb Z}.
\end{equation}
Here $f_{k, \ n}:=F_{k}(v(x),x)_{n}$. Evidently, the Fourier coefficients of 
the second derivatives are given by
$$
(-u_{k}'')_{n}=\sqrt{2\pi}\frac{n^{2}G_{k, \ n}f_{k, \ n}}{|n|-a_{k}}, \quad
a_{k}\geq 0, \quad 1\leq k\leq N_{1}, \quad n\in {\mathbb Z},
$$
$$
(-u_{k}'')_{n}=\sqrt{2\pi}\frac{n^{2}G_{k, \ n}f_{k, \ n}}{|n|+a_{k}}, \quad
a_{k}>0, \quad N_{1}+1\leq k\leq N_{2}, \quad n\in {\mathbb Z}.
$$
Let us use quantities (\ref{Pk1}), (\ref{Pk2}), (\ref{Pk3}), (\ref{Pk4}) and
(\ref{p}) to arrive at
$$
|u_{k, \ n}|\leq \sqrt{2\pi}P|f_{k, \ n}|, \quad |(-u_{k}'')_{n}|\leq
\sqrt{2\pi}P|f_{k, \ n}|, \quad 1\leq k\leq N_{2}, \quad n\in {\mathbb Z}.
$$
We easily derive the upper bound
\begin{equation}
\label{hc2iub}  
\|u\|_{H_{c}^{2}(I, \ {\mathbb R}^{N_{2}})}^{2}=\sum_{k=1}^{N_{2}}\Bigg\{\sum_{n=-\infty}^{\infty}
|u_{k, \ n}|^{2}+\sum_{n=-\infty}^{\infty}|(-u_{k}'')_{n}|^{2}\Bigg \}\leq
$$
$$
4\pi P^{2}\sum_{k=1}^{N_{2}}\sum_{n=-\infty}^{\infty}|f_{k, \ n}|^{2}=
4\pi P^{2}\sum_{k=1}^{N_{2}}\|F_{k}(v(x),x)\|_{L^{2}(I)}^{2}.
\end{equation}
Let us use inequality (\ref{ub1}) of Assumption 1.1. Hence, all
$F_{k}(v(x),x)\in L^{2}(I)$, such that the right side of (\ref{hc2iub}) is
finite.
This means that for an arbitrarily chosen vector function 
$v(x)\in H_{c}^{2}(I, \ {\mathbb R}^{N_{2}})$ 
there exists a unique $u(x)\in H_{c}^{2}(I, \ {\mathbb R}^{N_{2}})$
satisfying system (\ref{ae1}),  (\ref{ae2}) and its Fourier
coefficients are given by (\ref{f2}), (\ref{f21}), so that the map 
$\tau_{a}: H_{c}^{2}(I, \ {\mathbb R}^{N_{2}})\to H_{c}^{2}(I, \ {\mathbb R}^{N_{2}})$
is well defined.

\noindent
Obviously, orthogonality conditions (\ref{or4}) and (\ref{or5})
along with (\ref{f2}) imply that for $m+1\leq k\leq q$ the components
$u_{k}(x)$ are orthogonal to the Fourier harmonics 
$\displaystyle{\frac{e^{\pm in_{k}x}}{\sqrt{2\pi}}}$ in $L^{2}(I)$. Similarly, 
for $q+1\leq k\leq N_{1}$ the functions $u_{k}(x)$ are orthogonal to $1$ in
$L^{2}(I)$ since the corresponding Fourier coeffients are trivial.

\noindent
We consider arbitary vector functions 
$v^{(1),(2)}(x)\in H_{c}^{2}(I, \ {\mathbb R}^{N_{2}})$. Their
images under the map introduced above are
$u^{(1),(2)}:=\tau_{a}v^{(1),(2)}\in H_{c}^{2}(I, \ {\mathbb R}^{N_{2}})$. By means
of (\ref{ae1}), (\ref{ae2}) with $\Omega=I$, 
\begin{equation}
\label{ae11i}
\sqrt{-\frac{d^{2}}{dx^{2}}}u_{k}^{(1)}-a_{k}u_{k}^{(1)}=\int_{0}^{2\pi}G_{k}(x-y)F_{k}
(v_{1}^{(1)}(y),v_{2}^{(1)}(y),..., v_{N_{2}}^{(1)}(y), y)dy, \quad a_{k}\geq 0
\end{equation}
when $1\leq k\leq N_{1}$,
\begin{equation}
\label{ae21i}
\sqrt{-\frac{d^{2}}{dx^{2}}}u_{k}^{(1)}+a_{k}u_{k}^{(1)}=\int_{0}^{2\pi}G_{k}(x-y)
F_{k}(v_{1}^{(1)}(y),v_{2}^{(1)}(y),..., v_{N_{2}}^{(1)}(y), y)dy, \quad a_{k}>0
\end{equation}
when $N_{1}+1\leq k\leq N_{2}$. Analogously,
\begin{equation}
\label{ae31i}
\sqrt{-\frac{d^{2}}{dx^{2}}}u_{k}^{(2)}-a_{k}u_{k}^{(2)}=\int_{0}^{2\pi}G_{k}(x-y)F_{k}
(v_{1}^{(2)}(y),v_{2}^{(2)}(y),..., v_{N_{2}}^{(2)}(y), y)dy, \quad a_{k}\geq 0
\end{equation}
if $1\leq k\leq N_{1}$,
\begin{equation}
\label{ae41i}
\sqrt{-\frac{d^{2}}{dx^{2}}}u_{k}^{(2)}+a_{k}u_{k}^{(2)}=\int_{0}^{2\pi}G_{k}(x-y)
F_{k}(v_{1}^{(2)}(y),v_{2}^{(2)}(y),..., v_{N_{2}}^{(2)}(y), y)dy, \quad a_{k}>0
\end{equation}
if $N_{1}+1\leq k\leq N_{2}$.

\noindent
Let us apply the Fourier transform (\ref{fti}) to
both sides of systems (\ref{ae11i}), (\ref{ae21i}) and (\ref{ae31i}),
(\ref{ae41i}). This yields
\begin{equation}
\label{ukn11}  
u_{k, \ n}^{(1)}=\sqrt{2\pi}\frac{G_{k, \ n}f_{k, \ n}^{(1)}}{|n|-a_{k}}, \quad
a_{k}\geq 0, \quad 1\leq k\leq N_{1}, \quad n\in {\mathbb Z},
\end{equation}
\begin{equation}
\label{ukn21}  
u_{k, \ n}^{(1)}=\sqrt{2\pi}\frac{G_{k, \ n}f_{k, \ n}^{(1)}}{|n|+a_{k}}, \quad
a_{k}>0, \quad N_{1}+1\leq k\leq N_{2}, \quad n\in {\mathbb Z}
\end{equation}
and
\begin{equation}
\label{ukn31}  
u_{k, \ n}^{(2)}=\sqrt{2\pi}\frac{G_{k, \ n}f_{k, \ n}^{(2)}}{|n|-a_{k}}, \quad
a_{k}\geq 0, \quad 1\leq k\leq N_{1}, \quad n\in {\mathbb Z},
\end{equation}
\begin{equation}
\label{ukn41}  
u_{k, \ n}^{(2)}=\sqrt{2\pi}\frac{G_{k, \ n}f_{k, \ n}^{(2)}}{|n|+a_{k}}, \quad
a_{k}>0, \quad N_{1}+1\leq k\leq N_{2}, \quad n\in {\mathbb Z}.
\end{equation}
Here $f_{k, \ n}^{(1), (2)}$ stand for the Fourier images of
$F_{k}(v^{(1), (2)}(x), x)$ under transform (\ref{fti}).

\noindent
We recall definitions (\ref{Pk1}),
(\ref{Pk2}), (\ref{Pk3}), (\ref{Pk4}) and (\ref{p}). By virtue of 
(\ref{ukn11}), (\ref{ukn21}), (\ref{ukn31}) and (\ref{ukn41}), we have
for $1\leq k\leq N_{2}, \ n\in {\mathbb Z}$ 
$$
|u_{k, \ n}^{(1)}-u_{k, \ n}^{(2)}|\leq \sqrt{2\pi}P|f_{k, \ n}^{(1)}-f_{k, \ n}^{(2)}|,
$$
$$
|n^{2}u_{k, \ n}^{(1)}-n^{2}u_{k, \ n}^{(2)}|\leq \sqrt{2\pi}P
|f_{k, \ n}^{(1)}-f_{k, \ n}^{(2)}|.
$$
Thus,
$$
\|u^{(1)}-u^{(2)}\|_{H_{c}^{2}(I, \ {\mathbb R}^{N_{2}})}^{2}= \sum_{k=1}^{N_{2}}\Bigg \{
 \sum_{n=-\infty}^{\infty}
|u_{k, \ n}^{(1)}-u_{k,\ n}^{(2)}|^{2}+
\sum_{n=-\infty}^{\infty}|n^{2}(u_{k, \ n}^{(1)}-u_{k, \ n}^{(2)})|^{2}\Bigg \}\leq
$$
$$
4\pi P^{2}\sum_{k=1}^{N_{2}}\sum_{n=-\infty}^{\infty}|f_{k, \ n}^{(1)}-f_{k,\ n}^{(2)}|^{2}=
4\pi P^{2} \sum_{k=1}^{N_{2}}\|F_{k}(v^{(1)}(x),x)-F_{k}(v^{(2)}(x),x)\|_
{L^{2}(I)}^{2}.
$$
Clearly, by means of the Sobolev embedding theorem we have
$v^{(1),(2)}_{k}(x)\in H^{2}(I)\subset L^{\infty}(I)$ for $1\leq k\leq N_{2}$.

\noindent
Let us use inequality (\ref{lk1}). We arrive at
$$
\sum_{k=1}^{N_{2}}\|F_{k}(v^{(1)}(x),x)-F_{k}(v^{(2)}(x),x)\|_{L^{2}(I)}^{2}\leq L^{2}
\sum_{k=1}^{N_{2}}\|v_{k}^{(1)}(x)-v_{k}^{(2)}(x)\|_{L^{2}(I)}^{2},
$$
so that
\begin{equation}
\label{taui}  
\|\tau_{a}v^{(1)}-\tau_{a}v^{(2)}\|_{H_{c}^{2}(I, \ {\mathbb R}^{N_{2}})}\leq 2\sqrt{\pi}
PL\|v^{(1)}-v^{(2)}\|_{H_{c}^{2}(I, \ {\mathbb R}^{N_{2}})}.
\end{equation}
Note that the constant in the right side of estimate (\ref{taui}) is less than
one as assumed.
The Fixed Point Theorem gives us the 
existence and uniqueness of a vector function 
$v_{a}(x)\in H_{c}^{2}(I, \ {\mathbb R}^{N_{2}})$, which satisfies 
$\tau_{a}v_{a}=v_{a}$.
This is the only stationary solution of system 
(\ref{id1}), (\ref{id11}) in 
$H_{c}^{2}(I, \ {\mathbb R}^{N_{2}})$.

\noindent
Finally, we suppose that $v_{a}(x)$ vanishes identically
in the interval $I$. This will contradict to the stated assumption
that the Fourier coefficients
$G_{k, \ n}F_{k}(0,x)_{n}\neq 0$ for a certain $1\leq k\leq N_{2}$ and some
$n\in {\mathbb Z}$.  \hfill\lanbox


 \setcounter{equation}{0}

 \section{The System in the Layer Domain}

\bigskip

{\it Proof of Theorem 1.7.} Let us suppose that there exists 
$v(x)\in H^{2}(\Omega, \ {\mathbb R}^{N_{2}})$, which generates
$u^{(1),(2)}(x)\in H^{2}(\Omega, \ {\mathbb R}^{N_{2}})$ satisfying system 
(\ref{ae1}), (\ref{ae2}). Clearly, the difference of these vector
functions
$w(x):=u^{(1)}(x)-u^{(2)}(x)\in H^{2}(\Omega, \ {\mathbb R}^{N_{2}})$ will solve
the homogeneous system of equations 
\begin{equation}
\label{hl1}  
\sqrt{-\Delta} w_{k}=a_{k}w_{k}, \quad a_{k}\geq 0, \quad  1\leq k\leq N_{1},
\end{equation}
\begin{equation}
\label{hl2}  
\sqrt{-\Delta} w_{k}=-a_{k}w_{k}, \quad a_{k}>0, \quad  N_{1}+1\leq k\leq N_{2}.
\end{equation}
We apply the partial Fourier transform to (\ref{hl1}). This yields
$$
\sqrt{-\Delta_{\perp}+n^{2}}w_{k, \ n}(x_{\perp})=a_{k}w_{k, \ n}(x_{\perp}), \quad
1\leq k\leq N_{1}, \quad n\in {\mathbb Z}.
$$ 
Here 
$\displaystyle{w_{k, \ n}(x_{\perp}):=\frac{1}{\sqrt{2\pi}}\int_{0}^{2\pi}
w_{k}(x_{1},x_{\perp})e^{-inx_{1}}dx_{1}}$ and $\Delta_{\perp}$ stands for the
Laplacian with respect to the $x_{\perp}$ variable.
Evidently,
$$
\|w_{k}\|_{L^{2}(\Omega)}^{2}=\sum_{n=-\infty}^{\infty}\|w_{k, \ n}\|_
{L^{2}({\mathbb R}^{d})}^{2}, \quad 1\leq k\leq N_{1}.
$$ 
Hence,
$$
w_{k, \ n}(x_{\perp})\in L^{2}({\mathbb R}^{d}), \quad 1\leq k\leq N_{1}, \quad
n\in {\mathbb Z}.
$$
Obviously, each operator  $\sqrt{-\Delta_{\perp}+n^{2}}$ considered on
$L^{2}({\mathbb R}^{d})$ does not possess any nontrivial 
eigenfunctions. This means that $w_{k}(x)$ vanish in $\Omega$ for
$1\leq k\leq N_{1}$.

\noindent
The analogous assertion is valid for $w_{k}(x), \
N_{1}+1\leq k\leq N_{2}$ by means of (\ref{hl2}), because the operator
$\sqrt{-\Delta}: H^{2}(\Omega)\to L^{2}(\Omega)$ cannot have any negative
eigenvalues.

\noindent
We choose an arbitrary vector function 
$v(x)\in H^{2}(\Omega, \ {\mathbb R}^{N_{2}})$. Let us apply the Fourier transform
(\ref{ftl}) to both sides of system (\ref{ae1}), (\ref{ae2}) to arrive at
\begin{equation}
\label{f3}
\widehat{{u}_{k, \ n}}(p)=(2\pi)^{\frac{d+1}{2}}\frac{\widehat{{G}_{k, \ n}}(p)
\widehat{{f}_{k, \ n}}(p)}{\sqrt{p^{2}+n^{2}}-a_{k}}, \quad a_{k}\geq 0, \quad
1\leq k\leq N_{1}, \quad n\in {\mathbb Z}, \quad p\in
{\mathbb R}^{d},
\end{equation}
\begin{equation}
\label{f4}
\widehat{{u}_{k, \ n}}(p)=(2\pi)^{\frac{d+1}{2}}\frac{\widehat{{G}_{k, \ n}}(p)
\widehat{{f}_{k, \ n}}(p)}{\sqrt{p^{2}+n^{2}}+a_{k}}, \quad a_{k}>0, \quad
N_{1}+1\leq k\leq N_{2}, \quad n\in {\mathbb Z}, \quad p\in
{\mathbb R}^{d}.
\end{equation}
Here $\widehat{{f}_{k, \ n}}(p)$ denotes the Fourier image of $F_{k}(v(x),x)$
under transform (\ref{ftl}).

\noindent
Clearly, we have the upper bounds in terms of the expressions given
by (\ref{R1}), (\ref{R2}), (\ref{R3}), (\ref{R4}) and (\ref{R}), such that
$$
|\widehat{{u}_{k, \ n}}(p)|\leq (2\pi)^{\frac{d+1}{2}}R
|\widehat{{f}_{k, \ n}}(p)|, \quad 
|(p^{2}+n^{2})\widehat{{u}_{k, \ n}}(p)|\leq (2\pi)^{\frac{d+1}{2}}R
|\widehat{{f}_{k, \ n}}(p)|,
$$
where $1\leq k\leq N_{2}, \ n\in {\mathbb Z}, \ p\in {\mathbb R}^{d}$.
This yields
$$
\|u\|_{H^{2}(\Omega, \ {\mathbb R}^{N_{2}})}^{2}=\sum_{k=1}^{N_{2}}\Bigg\{
\sum_{n=-\infty}^{\infty}\int_{{\mathbb R}^{d}}
|\widehat{{u}_{k, \ n}}(p)|^{2}dp+
\sum_{n=-\infty}^{\infty}\int_{{\mathbb R}^{d}}
|(p^{2}+n^{2})\widehat{{u}_{k, \ n}}(p)|^{2}dp \Bigg\} \leq
$$
\begin{equation}
\label{uh2omubn}
2(2 \pi)^{d+1}R^{2}\sum_{k=1}^{N_{2}}\sum_{n=-\infty}^{\infty}\int_{{\mathbb R}^{d}}
|\widehat{{f}_{k, \ n}}(p)|^{2}dp=
2(2 \pi)^{d+1}R^{2}\sum_{k=1}^{N_{2}}
\|F_{k}(v(x),x)\|_{L^{2}(\Omega)}^{2}.
\end{equation}
Let us recall bound (\ref{ub1}) of Assumption 1.1. Hence, all
$F_{k}(v(x),x)\in L^{2}(\Omega)$, so that the right side
of (\ref{uh2omubn}) is finite.
Therefore, for any vector function $v(x)\in H^{2}(\Omega, \ {\mathbb R}^{N_{2}})$ 
there exists a unique $u(x)\in H^{2}(\Omega, \ {\mathbb R}^{N_{2}})$, which 
satisfies the system of equations (\ref{ae1}), (\ref{ae2}).
Its Fourier transform is given by (\ref{f3}), (\ref{f4}). This means that the
map
$t_{a}: H^{2}(\Omega, \ {\mathbb R}^{N_{2}})\to H^{2}(\Omega, \ {\mathbb R}^{N_{2}})$
is well defined.

\noindent
Let us choose arbitrarily two vector functions 
$v^{(1),(2)}\in H^{2}(\Omega, \ {\mathbb R}^{N_{2}})$.
Their images under the map discussed above are
$u^{(1),(2)}:=t_{a}v^{(1),(2)}\in H^{2}(\Omega, \ {\mathbb R}^{N_{2}})$. By virtue of
(\ref{ae1}), (\ref{ae2}), we have 
\begin{equation}
\label{ael1}
\sqrt{-\Delta} u_{k}^{(1)}-a_{k}u_{k}^{(1)}=\int_{\Omega}G_{k}(x-y)F_{k}
(v_{1}^{(1)}(y),v_{2}^{(1)}(y),..., v_{N_{2}}^{(1)}(y), y)dy, \quad a_{k}\geq 0
\end{equation}
for $1\leq k\leq N_{1}$,
\begin{equation}
\label{ael2}
\sqrt{-\Delta} u_{k}^{(1)}+a_{k}u_{k}^{(1)}=\int_{\Omega}G_{k}(x-y)
F_{k}(v_{1}^{(1)}(y),v_{2}^{(1)}(y),..., v_{N_{2}}^{(1)}(y), y)dy, \quad a_{k}>0
\end{equation}
for $N_{1}+1\leq k\leq N_{2}$. Similarly,
\begin{equation}
\label{ael3}
\sqrt{-\Delta} u_{k}^{(2)}-a_{k}u_{k}^{(2)}=\int_{\Omega}G_{k}(x-y)F_{k}
(v_{1}^{(2)}(y),v_{2}^{(2)}(y),..., v_{N_{2}}^{(2)}(y), y)dy, \quad a_{k}\geq 0
\end{equation}
when $1\leq k\leq N_{1}$,
\begin{equation}
\label{ael4}
\sqrt{-\Delta} u_{k}^{(2)}+a_{k}u_{k}^{(2)}=\int_{\Omega}G_{k}(x-y)
F_{k}(v_{1}^{(2)}(y),v_{2}^{(2)}(y),..., v_{N_{2}}^{(2)}(y), y)dy, \quad a_{k}>0
\end{equation}
when $N_{1}+1\leq k\leq N_{2}$.

\noindent
We apply the Fourier transform (\ref{ftl}) to
both sides of the systems of equations (\ref{ael1}), (\ref{ael2}) and
(\ref{ael3}), (\ref{ael4}). This gives us
\begin{equation}
\label{f31}
\widehat{{u}_{k, \ n}^{(1)}}(p)=(2\pi)^{\frac{d+1}{2}}\frac{\widehat{{G}_{k, \ n}}(p)
\widehat{{f}_{k, \ n}^{(1)}}(p)}{\sqrt{p^{2}+n^{2}}-a_{k}}, \quad a_{k}\geq 0, \quad
1\leq k\leq N_{1}, \quad n\in {\mathbb Z}, \quad p\in
{\mathbb R}^{d},
\end{equation}
\begin{equation}
\label{f41}
\widehat{{u}_{k, \ n}^{(1)}}(p)=(2\pi)^{\frac{d+1}{2}}\frac{\widehat{{G}_{k, \ n}}(p)
\widehat{{f}_{k, \ n}^{(1)}}(p)}{\sqrt{p^{2}+n^{2}}+a_{k}}, \quad a_{k}>0, \quad
N_{1}+1\leq k\leq N_{2}, \quad n\in {\mathbb Z}, \quad p\in
{\mathbb R}^{d}
\end{equation}
and
\begin{equation}
\label{f51}
\widehat{{u}_{k, \ n}^{(2)}}(p)=(2\pi)^{\frac{d+1}{2}}\frac{\widehat{{G}_{k, \ n}}(p)
\widehat{{f}_{k, \ n}^{(2)}}(p)}{\sqrt{p^{2}+n^{2}}-a_{k}}, \quad a_{k}\geq 0, \quad
1\leq k\leq N_{1}, \quad n\in {\mathbb Z}, \quad p\in
{\mathbb R}^{d},
\end{equation}
\begin{equation}
\label{f61}
\widehat{{u}_{k, \ n}^{(2)}}(p)=(2\pi)^{\frac{d+1}{2}}\frac{\widehat{{G}_{k, \ n}}(p)
\widehat{{f}_{k, \ n}^{(2)}}(p)}{\sqrt{p^{2}+n^{2}}+a_{k}}, \quad a_{k}>0, \quad
N_{1}+1\leq k\leq N_{2}, \quad n\in {\mathbb Z}, \quad p\in
{\mathbb R}^{d}.
\end{equation}
In the formulas above $\widehat{{f}_{k, \ n}^{(1), (2)}}(p)$ designate the Fourier
images of $F_{k}(v^{(1), (2)}(x), x)$ under transform (\ref{ftl}).

\noindent
Let us recall
(\ref{R1}), (\ref{R2}), (\ref{R3}), (\ref{R4}) and (\ref{R}). By means of
(\ref{f31}), (\ref{f41}), (\ref{f51}) and (\ref{f61}), we obtain
for $1\leq k\leq N_{2}, \ n\in {\mathbb Z}, \ p\in {\mathbb R}^{d}$ that
$$
\Bigg|\widehat{{u}_{k, \ n}^{(1)}}(p)-\widehat{{u}_{k, \ n}^{(2)}}(p)\Bigg|\leq
(2\pi)^{\frac{d+1}{2}}R\Bigg|\widehat{{f}_{k, \ n}^{(1)}}(p)-
\widehat{{f}_{k, \ n}^{(2)}}(p)\Bigg|,
$$
$$
\Bigg|(p^{2}+n^{2})\Bigg[\widehat{{u}_{k, \ n}^{(1)}}(p)-
\widehat{{u}_{k, \ n}^{(2)}}(p)\Bigg]\Bigg|\leq
(2\pi)^{\frac{d+1}{2}}R\Bigg|\widehat{{f}_{k, \ n}^{(1)}}(p)-
\widehat{{f}_{k, \ n}^{(2)}}(p)\Bigg|.
$$
Therefore,
$$
\|u^{(1)}-u^{(2)}\|_{H^{2}(\Omega, \ {\mathbb R}^{N_{2}})}^{2}=
\sum_{k=1}^{N_{2}}\Bigg\{\sum_{n=-\infty}^{\infty}\int_{{\mathbb R}^{d}}\Bigg\{
\Bigg|\widehat{{u}_{k, \ n}^{(1)}}(p)-\widehat{{u}_{k, \ n}^{(2)}}(p)\Bigg|^{2}dp
+
$$
$$
\sum_{n=-\infty}^{\infty}\int_{{\mathbb R}^{d}}
\Bigg|(p^{2}+n^{2})(\widehat{{u}_{k, \ n}^{(1)}}-\widehat{{u}_{k, \ n}^{(2)}}(p))\Bigg
|^{2}dp\Bigg\}\leq
$$
$$
2(2 \pi)^{d+1}{R}^{2}\sum_{k=1}^{N_{2}}\sum_{n=-\infty}^{\infty}
\int_{{\mathbb R}^{d}}\Bigg|\widehat{{f}_{k, \ n}^{(1)}}(p)-\widehat{{f}_{k, \ n}^{(2)}}(p)
\Bigg|^{2}dp=
$$
$$
2(2 \pi)^{d+1}{R}^{2}\sum_{k=1}^{N_{2}}\|F_{k}(v^{(1)}(x),x)-F_{k}(v^{(2)}(x),x)\|
_{L^{2}(\Omega)}^{2}.
$$
According to the Sobolev embedding theorem,  
$v^{(1),(2)}_{k}(x)\in H^{2}(\Omega)\subset L^{\infty}(\Omega)$ with
$1\leq k\leq N_{2}$. 
We recall bound (\ref{lk1}) to derive
$$
\sum_{k=1}^{N_{2}}\|F_{k}(v^{(1)}(x),x)-F_{k}(v^{(2)}(x),x)\|_{L^{2}(\Omega)}^{2}\leq
L^{2}\sum_{k=1}^{N_{2}}\|v_{k}^{(1)}(x)-v_{k}^{(2)}(x)\|_{L^{2}(\Omega)}^{2},
$$
such that
\begin{equation}
\label{fpl}  
\|t_{a}v^{(1)}-t_{a}v^{(2)}\|_{H^{2}(\Omega, \ {\mathbb R}^{N_{2}})}
\leq \sqrt{2}(2\pi)^{\frac{d+1}{2}}RL
\|v^{(1)}-v^{(2)}\|_{H^{2}(\Omega, \ {\mathbb R}^{N_{2}})}.
\end{equation}
The constant in the right side of inequality (\ref{fpl}) is less than one
as we assume.
The Fixed Point Theorem implies the existence and
uniqueness of a vector function
$v_{a}(x)\in H^{2}(\Omega, \ {\mathbb R}^{N_{2}})$, for which
$t_{a}v_{a}=v_{a}$ holds. This is the only stationary solution of the 
system of equations (\ref{id1}), (\ref{id11}) in
$H^{2}(\Omega, \ {\mathbb R}^{N_{2}})$.

\noindent
Let us suppose that $v_{a}(x)$ is trivial in $\Omega$. 
This will give us the contradiction to the given condition that 
there exist $1\leq k\leq N_{2}$ and $n\in \mathbb Z$, such that
$\hbox{supp}\widehat{{F_{k}(0,x)}_{n}}(p)\cap \hbox{supp}\widehat{{G}_{k, \ n}}(p)$ is a set 
of nonzero Lebesgue measure in ${\mathbb R}^{d}$.             \hfill\lanbox

\bigskip


\section*{Acknowledgement} The author is grateful to Israel Michael Sigal
for the partial support by the NSERC grant NA 7901.

\bigskip


\end{document}